\documentclass[11pt]{article}
\input{epsf.sty}
\usepackage{graphicx}
\usepackage{indentfirst}
\usepackage{cprotect}
\begin{document}
\title{\LARGE\bf {Dynamical behavior
	 of Pielou's difference system with exponential term}
  \footnote{ Corresponding author. E-mail address:mouyang@xmut.edu.cn }}
 \date{}
 \author{Ouyang Miao$^{1,2,*}$, Qianhong Zhang$^{3}$ \\
        {\small $^1$School of Mathematics, Southwest Jiaotong University, }\\
         {\small Chengdu, Sichuan 610000,China}\\
        {\small $^2$School of Mathematics and Statistics, Xiamen University of Technology, }\\
        {\small Xiamen, Fujian 361000, China}\\
        {\small $^3$School of Mathematics and Statistics, Guizhou University of Finance and Economics,}\\
        {\small  Guiyang, Guizhou 550025, China}}

\maketitle
\par\noindent

\small  {\bf Abstract } \par
  In this paper, we investigate a type of  Pielou's difference system   with exponential term
   $$
  y_{n+1}=\frac{az_n}{p+z_{n}}e^{-y_n},\; \;\;\;\;\;\\
  z_{n+1}=\frac{by_n}{q+y_{n}}e^{-z_n}.
  $$
  where the parameters $a,b,p,q, $ are positive real numbers and the initial
  values $y_0,z_0$ are arbitrary nonnegative real numbers. Using the mean value theorem and Lyapunov functional skills, we obtained some sufficient
  conditions which guarantee the boundedness and persistence of solution, and global asymptotic stability of the equilibriums. Moreover, two numerical
   examples are given to  elaborate on the results.



\section{Motivation}\par
 \setcounter{equation}{0}
   \par
 Discrete time single-species model is the most appropriate mathematical description of life histories of organism whose  reproduction occurs only once a year during a very short season. These models are widely  used in fisheries and many  organisms [1-5].

   In 1965, Pielou's equation, as a discrete analogue of the delay logistic equation, was proposed by Pielou [6-7].
   \begin{equation}\label{1.1}
    	x_{n+1}=\frac{\alpha x_{n}}{1+x_{n-1}} ,n=0,1,2,...,
    \end{equation}
    where $p$ is a positive real number.

    In fact, there are many conclusions on Pielou's equation that can be accounted briefly: If $\alpha \le 1$, then zero is a globally asymptotically stable equilibrium of Eq.(1.1) , if $\alpha > 1$, and $x_0\in(0,\infty)$, then $\bar{y}= \alpha-1$ is a positive globally asymptotically stable equilibrium of Eq(1.1) . Other properties, such as periodicity, oscillation of Pielou's equation, one can refer to [8-12].

    In 2001, Metwally, Grove and Ladas [13] investigated the global stability of a difference model in mathematical biology
    \begin{equation}\label{1.2}
    	x_{n+1}=\alpha+\beta x_{n-1} \mathrm{e}^{-x_{n}}, \quad n=0,1,2, \ldots,
    \end{equation}
where $\alpha$ is the immigration rate and $\beta $the population growth rate.

   In 2006, Ozturk [14] researched the convergence, the boundedness and periodic character of the positive solutions of the following difference equation
\begin{equation}
x_{n+1}=\frac{\alpha+\beta e^{-x_{n}}}{\gamma+x_{n-1}}, \quad n=0,1,2, \ldots,
\end{equation}
where the constants $\alpha, \beta, \gamma$, and the initial values $x_{-1}, x_0,$ are positive real numbers.

Based on the previous works, Papaschinopoulos et al. [15] studied the asymptotic behavior of the solutions of two difference equations systems with exponential form, respectively.
\begin{equation}\label{1.4}
\begin{array}{cc}
	  x_{n+1}=a+b y_{n-1} e^{-x_{n}}, \quad &
	     y_{n+1}=c+d x_{n-1} e^{-y_{n}} .   \\
	     \\
	  x_{n+1}=a+b y_{n-1} e^{-y_{n}}, \quad &
	      y_{n+1}=c+d x_{n-1} e^{-x_{n}} .
\end{array}
\end{equation}
where the constants $a, b, c, d$, and the initial values $
x_{-1}, x_0, y_{-1}, y_0$ are  positive real numbers.

Also, Papaschinopoulos et al.[16] studied the asymptotic behavior of the
solutions of difference systems of exponential form, respectively,
\begin{equation}\label{1.5}
	\begin{array}{ll}
		x_{n+1}=\frac{\alpha+\beta e^{-y_{n}}}{\gamma+y_{n-1}}, & y_{n+1}=\frac{\delta+\epsilon e^{-x_{n}}}{\zeta+x_{n-1}}. \\
		\\
		x_{n+1}=\frac{\alpha+\beta e^{-y_{n}}}{\gamma+x_{n-1}}, & y_{n+1}=\frac{\delta+\epsilon e^{-x_{n}}}{\zeta+y_{n-1}}. \\
		\\
		x_{n+1}=\frac{\alpha+\beta e^{-x_{n}}}{\gamma+y_{n-1}}, & y_{n+1}=\frac{\delta+\epsilon e^{-y_{n}}}{\zeta+x_{n-1}}.
	\end{array}
\end{equation}
where  $\alpha, \beta, \gamma, \delta, \epsilon, \zeta$, and the initial values  $x_{-1}, x_{0}, y_{-1}, y_{0}$  are positive constants.

Hereafter, many researchers have worked out  colorful homologous-series results. Readers can refer to [17-30].

In this paper, by incorporating  Pielou's equation  to  extend exponential type  difference system, we study the dynamical behavior of the solutions to the following system
\begin{equation}\label{1.6}
	y_{n+1}=\frac{az_n}{p+z_{n}}e^{-y_n},\; \;\;\;\;\;\\
	z_{n+1}=\frac{by_n}{q+y_{n}}e^{-z_n}.
\end{equation}
where the parameters $a,\;b,\;p,\;q\; $ are positive numbers,  and the initial values $x_0, y_0$  are arbitrary nonnegative real numbers.

 The main aim of our work is to study the boundedness and persistence of positive solutions of system (1.6), and its asymptotic behavior. Furthermore, using the iteration method of  nonlinear difference equations and Lyapunov function skills,  we derive some conditions so that system (1.6) has a unique positive equilibrium, and global asymptotically stability of the equilibrium.

\section{Persistence and boundedness}
\setcounter{equation}{0}
 \vskip2mm
\par\noindent

In order to establish the persistence and boundedness of system (1.6), we introduce the following definitions and notations.

 \vskip2mm
\par\noindent
{\bf Definition 2.1 } {\it Let $(y_n,z_n)$ be an solution of system (1.6), the boundedness and persistence of  $(y_n,z_n)$ is defined as follows:\\

(i) The sequences  $\{y_n\}, \{z_n\}$ of positive solution are said to be bounded, if there exist  positive real numbers $N_1, N_2,$
such that $\mbox{supp} y_{n}\subset\left( 0,N_1\right], \mbox{supp} z_{n} \subset \left( 0,N_2\right],  n=1,2, \cdots.$\\

(ii) The sequences  $\{y_n\}, \{z_n\}$ of positive solution  are said to be persistence, if there exist  positive real numbers $M_1,M_2 ,$
such that $\mbox{supp} y_{n}\subset \left[ M_1,\infty\right), \mbox{supp} z_{n} \subset \left[ M_2,\infty\right), n=1,2, \cdots.$\\

(iii) The sequences  $\{y_n\}, \{z_n\}$ of positive  numbers  are said to be bounded and persistence if there exist positive real numbers
 $M_1,N_1, M_2,N_2$ such that $\mbox{supp} y_{n} \subset[M_1, N_1], \mbox{supp} z_{n}\subset[M_2, N_2],$ $ n=1,2, \cdots$. }

 \vskip2mm
 \par\noindent
  {\bf Theorem 2.1 } {\it Consider system (1.6), if $p, q\in(0,+\infty),  a, b\in(0,1), y_0,z_0\in(0,+\infty)$, then the following statements are true.\\
    
    (i)  the positive solution  ${y_n}, {z_n}$ of system (1.6) is bounded.\\
    
    (ii) the positive solution ${y_n}, {z_n}$ of system (1.6)  is persistence.}
    \vskip2mm
  \par\noindent
  {\bf Proof.} (i) Let $\{y_n\}, \{z_n\} $ be the solution sequence of system (1.6), for
  any positive initial values, obviously, ${y_n}>0, {z_n}>0. $
  Therefore
    \begin{equation}\label{2.1}
    	e^{-y_n}\leq 1, \   \   e^{-z_n}\leq 1.
    \end{equation}
   From (1.6) and (2.1), we have
   \begin{equation}\label{2.2}
   	\left\{\begin{array}{l}
   		y_{n+1} \leq \frac{a z_{n}}{p+z_{n}}, \\
   		\\
   		z_{n+1} \leq \frac{b y_{n}}{q+y_{n}}.
   	\end{array}\right.
   \end{equation}
  Consider the system of difference equations
   \begin{equation}\label{2.3}
  	\left\{\begin{array}{l}
  		x_{n+1} = \frac{a w_{n}}{p+w_{n}}, \\
  		\\
  		w_{n+1} = \frac{b x_{n}}{q+x_{n}}.
  	\end{array}\right.
  \end{equation}
Let $(x_n, w_n) $ be a solution of (2.3) such that
\begin{equation}\label{2.4}
	x_{0}=y_{0} ,   \                \  w_{0}=z_{0}.
\end{equation}
Simplify (2.3),  we get
 \begin{equation}\label{2.4}
	\left\{\begin{array}{l}
	x_{n+1}=\frac{a b x_{n-1}}{pq+(p+b) x_{n-1}}, \\
	\\
	w_{n+1}=\frac{a b w_{n-1}}{pq+(q+a) w_{n-1}}.
\end{array}
\right.
\end{equation}
Set
$$X_{n}=\frac{1}{x_{n}},\   W_{n}=\frac{1}{w_{n}}.$$
(2.4) turns to
\begin{equation}\label{2.6}
	X_{0}=\frac{1}{y_{0}} ,   \                \  W_{0}=\frac{1}{z_{0}}.
\end{equation}
Then (2.5) turns into
\begin{equation}
		\left\{\begin{array}{l}
			X_{n+1}=\frac{p+b}{a b}+\frac{pq}{a b} X_{n-1},\\
	\\
	W_{n+1}=\frac{q+a}{a b}+\frac{pq}{a b} W_{n-1}.
\end{array}\right.\label{2.7}
\end{equation}
From (2.6) and (2.7), we obtain
\begin{equation}\label{2.8}
\left\{	\begin{array}{l}
		X_{n}=\lambda_{1}\left(\sqrt{\frac{pq}{a b}}\right)^{n}+\lambda_{2}\left(-\sqrt{\frac{pq}{a b}}\right)^{n}+\frac{p+b}{a b-pq}, \\
		\\
		W_{n}=\mu_{1}\left(\sqrt{\frac{pq}{a b}}\right)^{n}+\mu_{2}\left(-\sqrt{\frac{pq}{a b}}\right)^{n}+\frac{q+a}{a b-pq}.
	\end{array}\right.\label{2.8}
\end{equation}
Where $\lambda_{1}, \lambda_{2}, \mu_{1}, \mu_{2}$  are expressions of $X_0, W_0$.

From the relations between (2.2) and  (2.7), we know
\begin{equation}\label{2.9}
\left\{\begin{array}{l}
	y_{n+1} \leq  \frac{1}{\lambda_{1}\left(\sqrt{\frac{pq}{a b}}\right)^{n}+\lambda_{2}\left(-\sqrt{\frac{pq}{a b}}\right)^{n}+\frac{p+b}{a b-pq}},\\
	\\
	z_{n+1} \leq  \frac{1}{\mu_{1}\left(\sqrt{\frac{pq}{a b}}\right)^{n}+\mu_{2}\left(-\sqrt{\frac{pq}{a b}}\right)^{n}+\frac{q+a}{a b-pq}}.

	\end{array}\right.
\end{equation}
From Eq.(1.6), without lose of correctness, we can get  some loose upper bound to simplfy the expression,
\begin{equation}\label{2.10}
	\left\{\begin{array}{l}
		y_{n+1} \leq a,\\
		\\
		z_{n+1} \leq b.
	\end{array}\right.
\end{equation}
Therefore, we have the boundedness of the solution of Eq.(1.6) .

(ii) Consider the function
$$f(x)=\frac{a x}{p+ x},\  \  \   (\mbox{resp. }  g(x)= \frac{b x}{q +x}).$$
where  $x\in (0, +\infty).$
 Since $ f'(x)=\frac{ap}{(p+x)^2}>0, (\mbox{resp.}  g'(x)=\frac{bq}{(q+x)^2}>0,$  it's easy to know $f(x), g(x)$ are increasing functions. So, for $\forall (y_0,z_0) \in (0, +\infty)$
\begin{equation}\label{2.10}
	\left\{\begin{array}{l}
		y_{n+1} \ge \frac{az_0}{p
		+z_0} e^{-b},\\
	\\
        z_{n+1} \ge \frac{by_0}{q+y_0} e^{-a}.
\end{array}\right.
\end{equation}	
It shows the solutions of Eq.(1.6) is permanent.

 \vskip2mm
 \par\noindent
 {\bf Theorem 2.2. } {\it The positive solution ${y_n}, {z_n}$ of system (1.6) are boundedness and persistence.\\}

\section{Stability }
\setcounter{equation}{0}

In order to obtain the existence  and the attractivity of the
unique positive equilibrium of system (1.6), we give the following definitions, notations  and Lemmas.
\vskip2mm
 \par\noindent
{\bf Definition 3.1 } [14] {\it
 Let $I \subset R$ and let $ f,g: I \times I \rightarrow I $ be  continuous differential functions. Consider the systems of difference equations
\begin{equation}
y_{n+1}=f\left(y_{n}, z_{n}\right), z_{n+1}=g\left(y_{n}, z_{n}\right),\quad n=0,1,2, \cdots,
  \end{equation}
  where the initial conditions $ y_{0}, z_{0} \in I $.  $ (\bar{y},\bar{z}) $  is said to be an equilibrium of Eq. (3.1) if $$\bar{y}=f(\bar{y}, \bar{z}),  \bar{z}=g(\bar{y}, \bar{z}).$$}
\vskip2mm
 \par\noindent
{\bf Definition 3.2 } [14]
{ \it  (i) The equilibrium  $(\bar{y},\bar{z}) $ of Eq. (3.1) is called locally stable if for every  $\varepsilon>0 $, there exists  $\delta>0$  such that $ y_{0}, z_{0} \in I $ with  $\left|y_{0}-\bar{y}\right|+\left|z_{0}-\bar{z}\right|<\delta$ , then $ \left|y_{n}-\bar{y}\right| +  \left|z_{n}-\bar{z}\right| < \varepsilon$  for all $ n \ge 0	.$\\
(ii) The equilibrium  $(\bar{y},\bar{z})$ of Eq. (3.1) is called locally asymptotically stable if it is locally stable, and if there exists $ \gamma>0  $ such that $ y_{0}, z_{0} \in I $  with  $\left|y_{0}-\bar{y}\right|+\left|z_{0}-\bar{z}\right|<\gamma $, then $ \lim _{n \rightarrow  \infty} y_{n}=\bar{y},  \lim _{n \rightarrow \infty} z_{n}=\bar{z}.$ \\
(iii) The equilibrium $(\bar{y},\bar{z})$ of Eq. (3.1) is called a global attractor if for every $ y_{0}, z_{0} \in I$  we have  $\lim _{n \rightarrow  \infty} y_{n}=\bar{y},  \lim _{n \rightarrow \infty} z_{n}=\bar{z}$.\\
(iv) The equilibrium  $(\bar{y},\bar{z})$ of Eq. (3.1) is called globally asymptotically stable if it is locally stable and a global attractor.\\
(v) The equilibrium $(\bar{y},\bar{z})$ of Eq. (3.1) is unstable if it is not stable.
}
\\
	

\vskip2mm
 \par\noindent
{\bf Theorem 3.1.} \label{3.3} {\it  (i) System (1.6) always has a zero equilibrium.\\
(ii) If \begin{equation}
	\frac{ab}{pq} > 1,
\end{equation}
then system (1.6) has a unique positive equilibrium.}
\vskip2mm
 \par\noindent
{\bf Proof}
(i)  It's obvious. So we omit.\\
 (ii) Let  $(\bar{y}, \bar{z} )$ be the positive equilibrium of system (1.6), i.e.
 \begin{equation}
 	\bar{y}=\frac{a \bar{z}}{p+\bar{z}} e^{-\bar{y}}, \quad \bar{z}=\frac{b \bar{y}}{q+\bar{y}} e^{-\bar{z}}.
 \end{equation}
From (2.11), we know
\begin{equation}
0 \le	y_0 e^{-b} <  \bar{y}	<	 a,\        \
0 \le z_0 e^{-a}  	<  \bar{z}	<	 b.
\end{equation}
We consider the homologous algebraic system of system (1.6),
\begin{equation}
y=\frac{a z}{p+z} e^{-y}, \quad z=\frac{b y}{q+y} e^{-z}.	\end{equation}
From (3.5), one has that
\begin{equation}
	z=\frac{p y e^y}{a-y e^y} ,\   \  y=\frac{q z e^z}{b-z e^z}.
\end{equation}
Substituting (3.6) into  system (3.5), we get
\begin{equation}
	e^{\frac{p y e^y}{a-y e^y}}= \frac{b y}{q + y}\frac{a-ye^y	}{p y e^y}, \           \
    e^{\frac{q z e^z}{b-z e^z}}= \frac{a z}{p + z}\frac{b-ze^z}{q z e^z}.
\end{equation}
\\
Set
\begin{equation}
	G(y)= e^{\frac{p y e^y}{a-y e^y}}+ \frac{b }{q + y}\frac{ye^y-a
	}{p e^y},
\end{equation}
Let $\epsilon_1 \rightarrow 0^+$,
\begin{equation}
	G(\epsilon_1)= 1- \frac{ab }{pq },
\end{equation}
Since condition (3.2) is satisfied, we have
\begin{equation}
 G(\epsilon_1)<0.
\end{equation}
On the other hand, from (3.8),
\begin{equation}
G(a)= e^{\frac{p a e^y}{a-a e^a}}+ \frac{b }{q + a}\frac{ae^a-a}{p e^a} > 0.
\end{equation}
and
\begin{equation}
	G'(y)=e^{\frac{p y e^{y}}{a-y}} \cdot\frac{a p(y+1) e^{y}}{\left(a-y e^{y}\right)^{2}}+\frac{e^{y}+a}{p e^{y}} \cdot \frac{b}{q+y}+\frac{a-y e^{y}}{p e^{y}} \frac{b}{(q+y)^{2}}>0.
\end{equation}
(3.12) implies $G(y)$ is an increasing function for $y \in (\epsilon_1,a)$.\\
Similarly, set
	\begin{equation}
		H(z)= e^{\frac{q z e^z}{b-z e^z}}+ \frac{a }{p + z}\frac{ze^z-b
		}{q e^z},
	\end{equation}
Let $\epsilon_2 \rightarrow 0^+$,
\begin{equation}
 H(\epsilon_2)= 1- \frac{ab}{pq },
\end{equation}
Noting condition (3.2), we have
\begin{equation}
	H(\epsilon_2) < 0.
\end{equation}
 From (3.13), we have
 \begin{equation}
	H(b)= e^{\frac{q b e^z}{b-b e^b}}+ \frac{a }{p + b}\frac{be^b-b}{q e^b} > 0,
\end{equation}
and
\begin{equation}
	H'(z)=e^{\frac{q z e^{z}}{b-z}} \cdot\frac{b q(z+1) e^{z}}{\left(b-z e^{z}\right)^{2}}+\frac{e^{z}+b}{q e^{z}} \cdot \frac{a}{p+z}+\frac{b-z e^{z}}{q e^{z}} \frac{a}{(p+z)^{2}}>0.
\end{equation}
(3.17) shows $H(z)$  is increasing function for $z \in (\epsilon_2,b)$.\\

 Therefore, we have that $(\bar{y},\bar{z}) \in (\epsilon_1, a) \times  (\epsilon_2, b)$ is the unique positive equilibrium of (1.6).\\

To have a more accurate estimation of $(\bar{y},\bar{z})$, we give a more narrow range by (2.9). Here, noting  the lower upper bounds as $y^*, z^*$ ,
\begin{equation}
		y_{n+1} < \frac{a b-pq}{p+b}=y^*,\     \
			    z_{n+1} < \frac{a b-pq}{q+a}=z^*.
	\end{equation}
Set $(y_*,z_*)$ be higher bound of $(\bar{y},\bar{z})$ , i. e.,
\begin{equation}
	y_* =\frac{q z_* e^{z_*}}{b-z_* e^{z_*}}	< y_{n+1} ,\     \
	z_* = \frac{p y_* e^{y_*}}{a-y_* e^{y_*}}  < z_{n+1}.
	\end{equation}
In fact, considering function $g(x)=\frac{c x e^x}{d-x e^x}$, since $g'(x)=\frac{cd(x+1)e^x}{(d-x e^x)^2}>0$ for $   x \in(0, \infty)$.
That means  $y$ (resp. $z$) in (3.6) will attain its minimums when $z$ (resp. $y$) come to the minimum.\\

Now, let's find the expressions of $y_*, z_*$. \\
 Noting (\ref{1.6}), since $y_{n+1}$ increases responding $z_n$,  decrease responding $y_n$,  we have
 \begin{equation}
	 y_{n}\geq \frac{a z_*}{p+z_*} e^{-y^*}  , \   \
	 z_{n}\geq \frac{b y_*}{q+y_*} e^{-z^*}.
\end{equation} \\
By  the definition of equilibrium,  noting (3.18), there exists $N>0,$ for  $n>N,$  $y_n= y_*, z_n=z_*$,  i.e.,
\begin{equation}
	 y_* = \frac{a z_*}{p+z_*} e^{-y^*}  , \   \
	 z_* = \frac{b y_*}{q+y_*} e^{-z^*}.
\end{equation}
Manipulating (3.21), we have
 \begin{equation}
		y_{*}=\frac{a b e^{-\left(y^{*}+z^*\right)}-p q}{b e^{-z^{*}}+p} ,\        \
		z_{*}=\frac{a b e^{-\left(y^{*}+z^{*}\right)}-p q}{a e^{-y^{*}}+q}.
\end{equation}
The unique positive equilibrium $(\bar{y},\bar{z})$ of system (\ref{1.6}) 
\begin{equation} 
	y_*< \bar{y}< y^*,\      \   z_* < \bar{z} <  z^*.
\end{equation}
The proof is completed.\\
\vskip2mm
 \par\noindent
{\bf Theorem 3.2.}\label{3.2}
 {\it
(i) If \begin{equation}
	\frac{a}{p} < 1, \  \  \frac{b}{q} < 1.
\end{equation}$(0,0)$ is locally asymptotically stable.\\
\\
(ii) If
\begin{equation}
	e^{y_*} >  a,    \    \ e^{z_*} >  b.
\end{equation}
$(\bar{y},\bar{z} )$ is locally asymptotically stable.}
\vskip2mm
 \par\noindent
{\bf Proof.}  (i) From system (\ref{1.6}), the linearized equation of system (1.6) at $(0,0)$ is
 \begin{equation}
	\phi_{n+1}=T ~\phi_{n}.
\end{equation}
where
\begin{equation}
\begin{array}{ll}
\phi_{n}=\left[
\begin{array}{c}
	y_{n} \\ z_{n}
  \end{array}
  \right],
&
	T=\left[
   \begin{array}{cc}
		0 & \frac{a}{p} \\
		\frac{b}{q} & 0
  \end{array}
   \right].
   \end{array}
\end{equation}
Under condition (3.24), by Theorem 1.3.7  [25],  the equilibrium point (0,0) of system (1.6) is locally asymptotically  stable.\\

(ii) The  linearized system of system (\ref{1.6}) at $(\bar{y},\bar{z})$ is
\begin{equation}\label{3.27}
	\phi_{n+1}=T_{(\bar{y},\bar{z})} \phi_{n}=
	\left[\begin{array}{cc}
		-\frac{a  \bar{z} e^{- \bar{y} }}{p+ \bar{z} }  & \frac{a pe^{- \bar{y} } }{\left(p+ \bar{z} \right)^{2}}  \\
		\\
		\frac{b q e^{- \bar{z} }}{\left(q+ \bar{y} \right)^{2}} & -\frac{b  \bar{y} e^{- \bar{z} }}{q+ \bar{y} }  \end{array}\right] ~\phi_{n}.
	\end{equation}
where $$
    T_{(\bar{y},\bar{z})}=
   \left[
   \begin{array}{cc}
    A &B\\
	C &D
  \end{array}
  \right].$$
In fact,
\begin{equation}
\begin{array}{c}	
	|A|+|B|=\left|-\frac{a  \bar{z} e^{- \bar{y} }}{p+ \bar{z} }\right| +
	\left| \frac{a pe^{- \bar{y} } }{(p+ \bar{z} )^2}  \right|  =
	\frac{a  \bar{z} e^{- \bar{y} }(p+ \bar{z} )+a pe^{- \bar{y} } }{(p+ \bar{z} )^2}  <
	a e^{- \bar{y} }<  a e^{-y_*}, \\
	\\
	|C|+|D|=\left| \frac{b qe^{- \bar{z} } }{(q+ \bar{y} )^2}  \right| +
	\left|-\frac{b  \bar{y} e^{- \bar{z} }}{q+ \bar{y} }\right|  =
	\frac{b  \bar{y} e^{- \bar{z} }(q+ \bar{y} )+b qe^{- \bar{z} } }{(q+ \bar{y} )^2}  <
	b e^{- \bar{z} }<  b e^{-z_*}.
\end{array}
\end{equation}
 Noting condition (3.25),  the inequations (3.29)  implies 
 \begin{equation}
 		|A|+|B| < 1,\    \  |C|+|D| < 1.
 \end{equation}
 The eigenvalue equation of (3.28) is
 \begin{equation}
 	\lambda^2 +\left( A+D \right)   	\lambda + \left(  AD-BC \right) =0.
 \end{equation}
 From (3.30), one has
 \begin{equation}
 	|A+D| < 1+ (A D-B C) < 2.
 \end{equation}
 by Theorem 1.1.1 [26],  the eigenvalues of Eq.(3.31) lie  inside the unit disk.
 So  the positive equilibrium $(\bar{y},\bar{z}) $ of system (1.6) is locally asymptotically stable.\\
 
 The proof is completed.

\vskip2mm
 \par\noindent
{\bf Theorem 3.3. }\label{3.3} {\it
(i) If \begin{equation}
	a b < \min (p,q),
\end{equation}
 the zero equilibrium is a global attractor.\\
(ii) If
\begin{equation}
	a z^* \le \bar{y}\left(p+z^*\right) e^{y_*},\  \ ~~~~~~
	b y^* \le \bar{z}\left(q+y^*\right) e^{z_*}.
\end{equation}
the equilibrium $(\bar{y},\bar{z})$ is a global attractor.
}
\\
{\bf Proof.}  Consider the following discrete time analog of the Lyapunov
function, mentioned in [31].
$$
V_{n}(\bar{y},\bar{z})=\left[\left(y_{n}-\bar{y}\right)-1-\ln \left(y_{n}-\bar{y}\right)\right]+\left[\left(z_{n}-\bar{z}\right)-1-\ln \left(z_{n}-\bar{z}\right)\right].$$
The nonnegativity of $V_n$ follows from: 
$$x-1-\ln x \geq 0, \quad \forall x>0.$$\\
Herefore, we have$$-\ln \left(\frac{y_{n+1}-\bar{y}}{y_{n}-\bar{y}}\right) \leq-\frac{y_{n+1}-y_{n}}{y_{n+1}-\bar{y}},\    \  -\ln \left(\frac{z_{n+1}-\bar{z}}{z_{n}-\bar{z}}\right) \leq-\frac{z_{n+1}-z_{n}}{z_{n+1}-\bar{z}}.$$
Assume
\begin{equation}
\begin{array}{ll}
V_{n+1}(\bar{y},\bar{z})-V_{n}(\bar{y},\bar{z})=&
\left[\left(y_{n+1}-\bar{y}\right)-1-\ln \left(y_{n+1}-\bar{y}\right)\right]-\left[\left(z_{n+1}-\bar{z}\right)-1-\ln \left(z_{n+1}-\bar{z}\right)\right]\\
\\
 &-  \left[\left(y_{n}-\bar{y}\right)-1-\ln \left(y_{n}-\bar{y}\right)\right]+\left[\left(z_{n}-\bar{z}\right)-1-\ln \left(z_{n}-\bar{z}\right)\right]\\
\\
&\le \left(y_{n+1}-y_{n}\right)\left(1-\frac{1}{y_{n+1}-\bar{y}}\right)+\left(z_{n+1}-z_{n}\right)\left(1-\frac{1}{z_{n+1}-\bar{z}}\right)\\
\\
& \le  (a-\epsilon_1)(\frac{y_{n+1}-\bar{y}-1}{y_{n+1}-\bar{y}}) + (b-\epsilon_2)(\frac{z_{n+1}-\bar{z}-1}{z_{n+1}-\bar{z}}).
\end{array}	
\end{equation}
Consider the zero equilibrium  (0,0),
\begin{equation}
	\begin{array}{ll}
		V_{n+1}(0,0)-V_{n}(0,0)    & \le
		(a-\epsilon_1)(\frac{y_{n+1}-1}{y_{n+1}}) + (b-\epsilon_2)(\frac{z_{n+1}-1}{z_{n+1}}) \\
		\\
		&  =
		(a-\epsilon_1)(\frac{az_n -(p+z_{n})e^{y_n}}{az_n}) + (b-\epsilon_2)(\frac{by_n -(q+y_{n})e^{z_n}}{by_n})  \\
		\\
		& \le
		(a-\epsilon_1)(\frac{ab -(p+\epsilon_2)e^{\epsilon_1}}{a\epsilon_2}) + (b-\epsilon_2)(\frac{ba -(q+\epsilon_1)e^{\epsilon_2}}{b\epsilon_1}).	
\end{array}	
\end{equation}
Set $\epsilon=\min (\epsilon_1,\epsilon_2),$  (3.36) becomes
\begin{equation}
		V_{n+1}(0,0)-V_{n}(0,0) \le
		(a-\epsilon)(\frac{ab -(p+\epsilon)e^{\epsilon}}{a\epsilon}) + (b-\epsilon)(\frac{ab -(q+\epsilon)e^{\epsilon}}{b\epsilon}).
\end{equation}
 Let $ \epsilon \rightarrow 0^+$, since (3.33) holds true, we have
\begin{equation}
	V_{n+1}(0,0)-V_{n}(0,0) \le 0,
\end{equation}
for all $n \geq  0$.

For the fact $V_n$ is  non-increasing and non-negative, we know that $\lim_{n\rightarrow \infty} V_n \geq 0. $ 
  Hence,  $$\lim_{n\rightarrow \infty} (V_{n+1}-V_{n}) = 0. $$ It follows that $\lim_{n\rightarrow \infty} y_n= 0, \lim_{n\rightarrow \infty} z_n =0$.
  Furthermore, $V_n \le  V_0$, for all $n \geq 0$, which shows
 that $(0,0) \in [\epsilon_1, a] \times [\epsilon_2, b]$ is a global attractor. 

(ii )  Next, we will prove the globally asymptotically stability of the unique positive equilibrium $(\bar{y},\bar{z})$.\\
We consider the following discrete time analog of the Lyapunov function
$$
W_{n}(\bar{y},\bar{z})=\bar{y}\left(\frac{y_{n}}{\bar{y}}-1-\ln \frac{y_{n}}{\bar{y}}\right)+\bar{z}\left(\frac{z_{n}}{\bar{z}}-1-\ln \frac{z_{n}}{\bar{z}}\right).$$
\\
The nonnegativity of $W_n$ follows from the following inequality:
$$x-1-\ln x \geq 0, \quad \forall x>0.$$
Furthermore, 
$$-\ln \left(\frac{y_{n+1}}{y_{n}}\right) \leq-\frac{y_{n+1}-y_{n}}{y_{n+1}},\    \  -\ln \left(\frac{z_{n+1}}{z_{n}}\right) \leq-\frac{z_{n+1}-z_{n}}{z_{n+1}}.$$
Since (3.34) holds true,  it follows that
\begin{equation}\
	\begin{array}{ll}
		W_{n+1}-W_{n}=&\bar{y}\left(\frac{y_{n+1}}{\bar{y}}-1-\ln\frac{y_{n+1}}{\bar{y}}\right)+\bar{z}\left(\frac{z_{n+1}}{\bar{z}}-1-\ln \frac{z_{n+1}}{\bar{z}}\right)\\
		\\
		&-\bar{y}\left(\frac{y_{n}}{\bar{y}}-1-\ln \frac{y_{n}}{\bar{y}}\right)-\bar{z}\left(\frac{z_{n}}{\bar{z}}-1-\ln \frac{z_{n}}{\bar{z}}\right)  \\
		\\
		&\leq \left(y_{n+1}-y_{n}\right)\left(1-\frac{\bar{y}}{y_{n+1}}\right)+\left(z_{n+1}-z_{n}\right)\left(1-\frac{\bar{z}}{z_{n+1}}\right)\\
		\\
		& =  \left(y_{n+1}-y_{n}\right) \frac{a z_{n}-\bar{y}\left(p+z_{n}\right) e^{y_{n}}}{a z_{n}} +\left(z_{n+1}-z_{n}\right) \frac{b y_{n}-\bar{z}\left(q+y_{n}\right) e^{z_{n}}}{b y_{n}}\\
		\\
		&\le \left(y^*-y_*\right) \frac{a z^*-\bar{y}\left(p+z^*\right) e^{y_*}}{a z^*} +         \left(z^*-z_*\right) \frac{b y^*-\bar{z}\left(q+y^* \right) e^{z_*}}{b y^*}\\
		\\
		& \le 0.
	\end{array}
\end{equation}
In fact, $f(x)=\frac{p+x}{ax}$ is a non-increase function , $-\frac{p+x}{ax} \le -\frac{p+x^*}{ax^*},$ for $\forall  x\in (x_*,x^*) $. \\
As $W_n$ being  non-increasing and non-negative,  it follows that $\lim_{n\rightarrow \infty} W_n \geq 0. $
 Hence,  $$\lim_{n\rightarrow \infty} W_{n+1}-W_{n} = 0.$$
 It follows that $$\lim_{n\rightarrow \infty} y_n= \bar{y}, \lim_{n\rightarrow \infty} z_n =\bar{z}.$$
  Furthermore, $W_n \le  W_0$ for all $n \geq 0$, which shows
that $(\bar{y},\bar{z}) \in [y_*, y^*] \times [z_*, z^*]$ is a global attractor. \\

\vskip2mm
\par\noindent
{\bf Theorem 3.4. }\label{3.4} {\it
	(i) If (3.33), the zero equilibrium is global asympotically stable.\\
	(ii) If (3.34), the positive equilibrium  $(\bar{y},\bar{z})$ is global asympotically stable.\\
}
\\
{\bf Proof.}
 By Theorem 3.2, Theorem 3.3, the  equilibrium  $(0,0)$ and $(\bar{y}, \bar{z})$ of system (1.6) is globally asymptotically stable.

\section{Numerical examples}
\setcounter{equation}{0}

{\bf Example 4.1}  Consider the parameters  $a=0.8, b=0.9 , p=0.6, q=0.5,$ and initial values    $y(1)_0=0.35,z(1)_0=0.26,$ the  Pielou's system with exponential terms (1.6) is as follows,  \\
\begin{equation}\label{key}
y_{n+1}=\frac{0.8z_n}{0.6+z_{n}}e^{-y_n},\; \;\;\;\;\;\\
z_{n+1}=\frac{0.9y_n}{0.5+y_{n}}e^{-z_n}.
\end{equation}\\
The solution of Eq.(4.1) $(y_n,z_n) \in [\frac{az_0}{p+z_0}e^{-a},a]\times [\frac{by_0}{q+y_0}e^{-b},b]= [0.1087, 0.8]\times[ 0.1507,0.9]$.\\
For the initial value $y(2)_0=0.05,z(2)_0=0.02$ ,  the solution   $(y_n,z_n) \in[0.0116, 0.8]\times[ 0.0333,0.9]$.\\
\begin{figure} [htbp]
	\begin{minipage}[t]{0.5\linewidth}
		\centering
		\includegraphics[width=3 in]{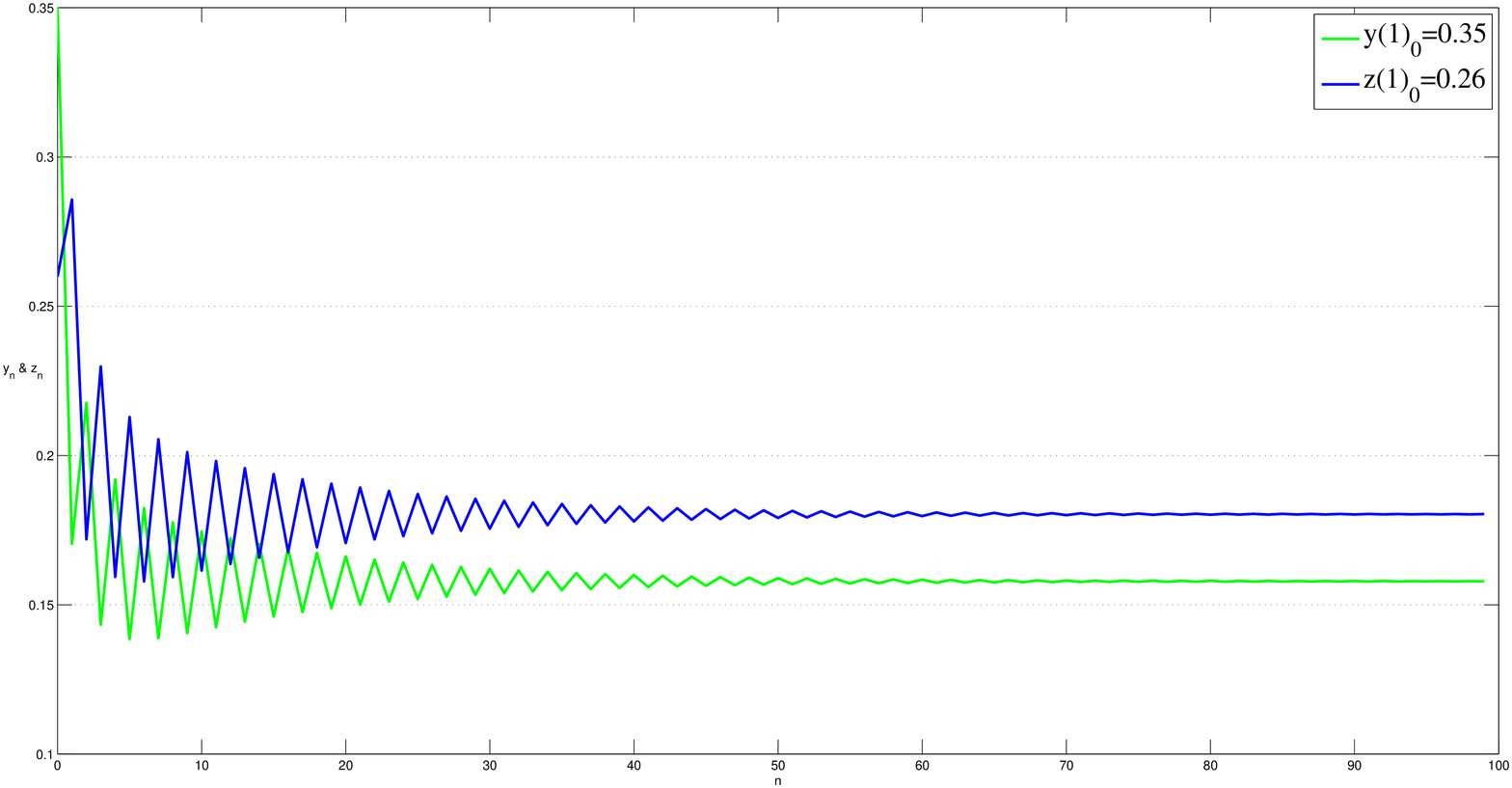}\\
		\centering{(a)  $y(1)_0=0.35,z(1)_0=0.26.$}
	\end{minipage}%
	\begin{minipage}[t]{0.5\linewidth}
		\centering
		\includegraphics[width=3 in]{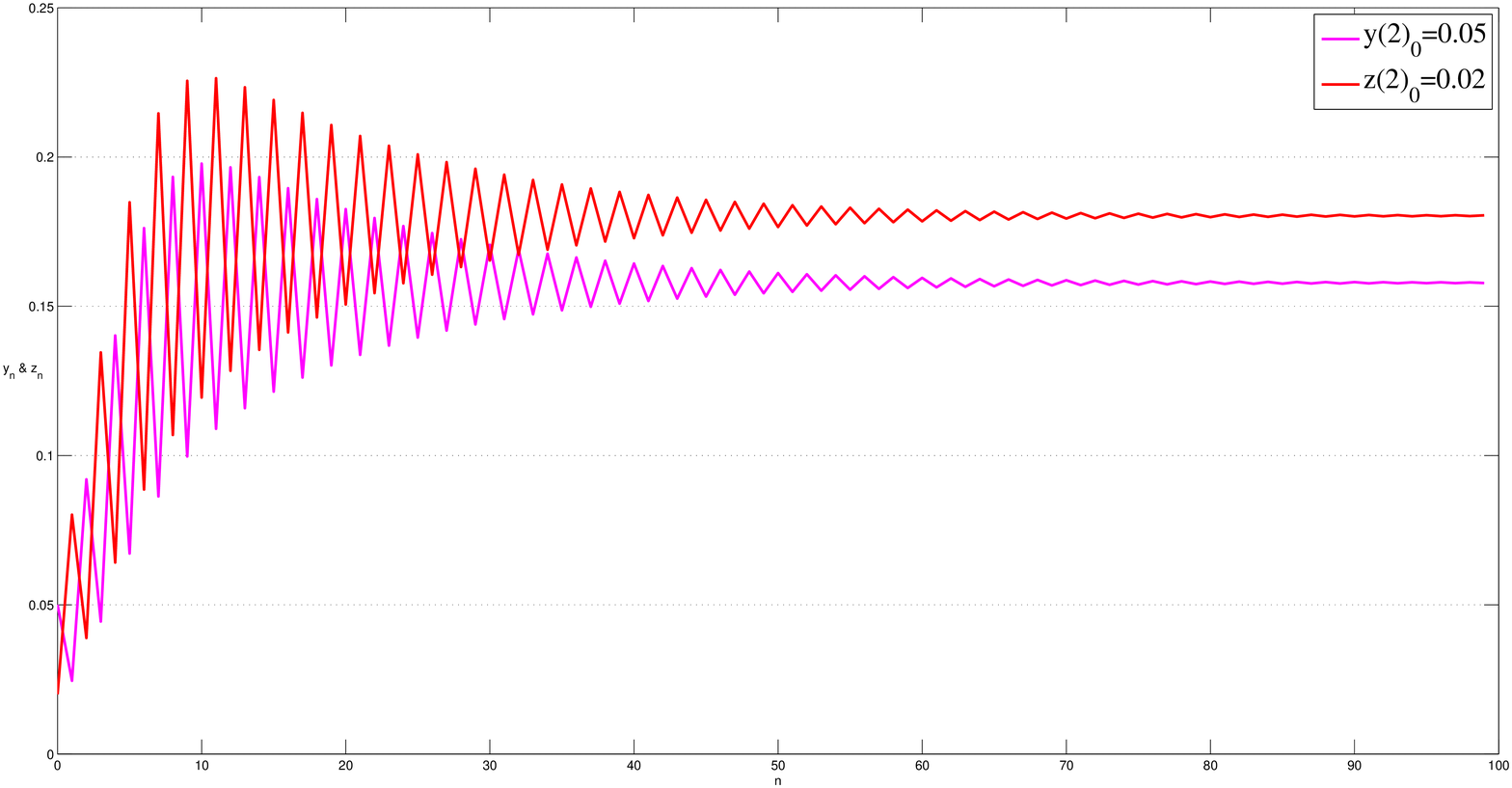}\\
		\centering{(b)$y(2)_0=0.05,z(2)_0=0.02.$}
	\end{minipage}
	\caption{The positive equilibrium dynamics of system (4.1).}
\end{figure}

It is obvious that conditions (3.2),(3.25) and (3.34) are satisfied. For the positive initial value wherever, Eq.(4.1) has the unique positive equilibrium   $(\bar{y},\bar{z}) \in [y_*,y^*]\times[z_*,z^*]$,(see Fig.1-Fig.2).\\ where 
$$
	y_{*}=\frac{a b e^{-\left(y^{*}+z^*\right)}-p q}{b e^{-z^{*}}+p} =0.0751,\     \ y^*= \frac{p+b}{a b-p^{2}}= 0.2800.
	\\
$$
$$
	z_{*}=\frac{a b e^{-\left(y^{*}+z^{*}\right)}-p q}{a e^{-y^{*}}+q} =0.0850,\     \ z^*= \frac{q+a}{a b-q^{2}}= 0.3231.
$$

\begin{figure}[htbp]
	\centering
	\includegraphics[totalheight=3in]{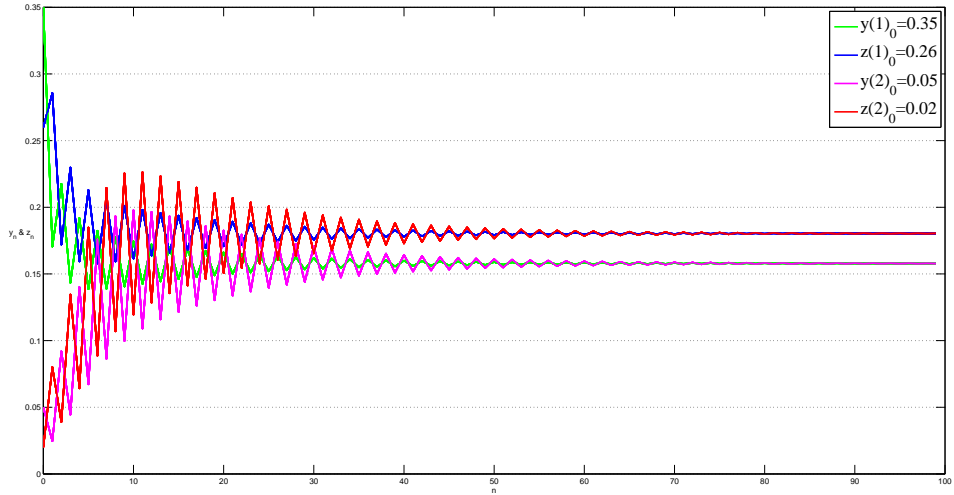}
	\caption{The comparation of solutions of system (4.1) from $(y(1)_0,z(1)_0)$ and $(y(2)_0,z(2)_0)$. } \label{fig:2}
\end{figure}

{\bf Example 4.2} Consider the parameters  $a=0.6, b=0.5 , p=0.8, q=0.9,$ the  Pielou's system with exponential terms (1.6) is as follows,  \\
\begin{equation}\label{key}
	y_{n+1}=\frac{0.6z_n}{0.8+z_{n}}e^{-y_n},\; \;\;\;\;\;\\
	z_{n+1}=\frac{0.5y_n}{0.9+y_{n}}e^{-z_n}.
\end{equation}\\
For the positive initial value wherever, the solution of Eq.(4.2) $(y_n,z_n) \in  [0, 0.6]\times[ 0,0.5]$.\\

The parameters meet  that conditions (3.24) and (3.33) . It is obvious that $(0,0)$ is the globally asymptotic stable. (see Fig.3-Fig.4.). 
\begin{figure} [htbp] 
	\begin{minipage}[t]{0.5\linewidth}
		\centering
		\includegraphics[width=3 in]{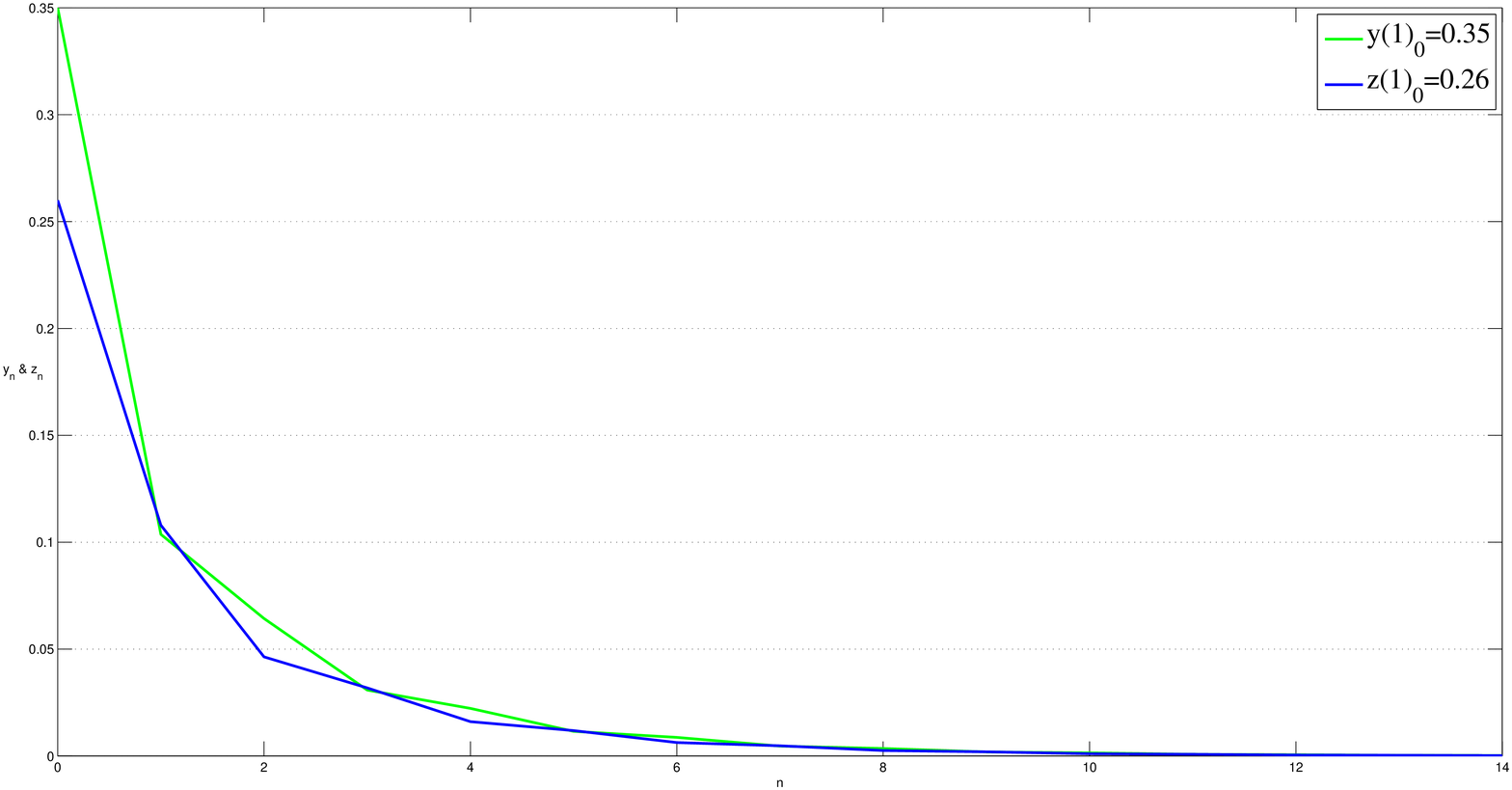}\\
		\centering{(a)  $y(1)_0=0.35,z(1)_0=0.26.$}
	\end{minipage}%
	\begin{minipage}[t]{0.5\linewidth}
		\centering
		\includegraphics[width=3 in]{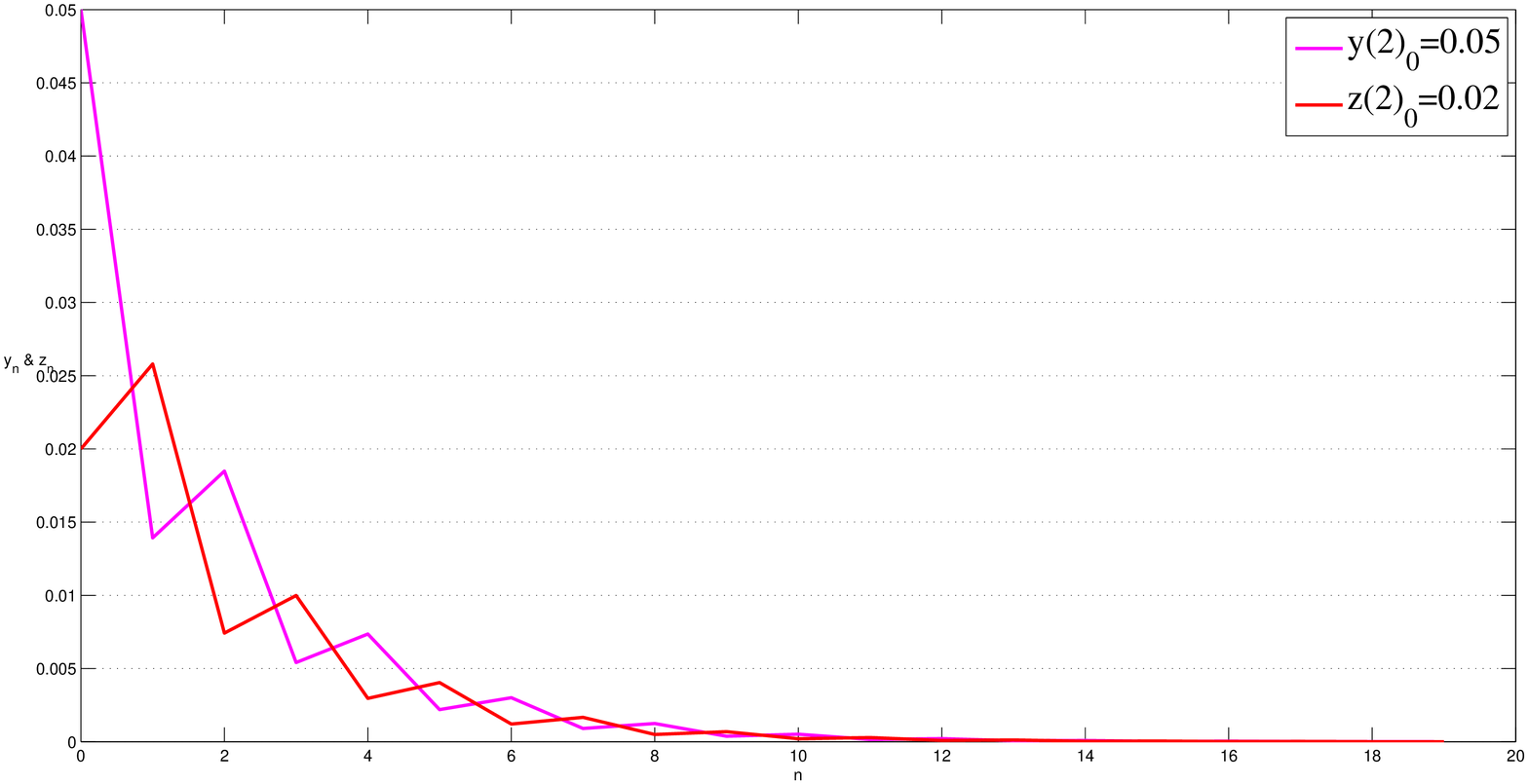}\\
		\centering{(b)  $y(2)_0=0.05,z(2)_0=0.02.$}
	\end{minipage}
	\caption{The zero equilibrium dynamics of system (4.2).}
\end{figure}

\begin{figure}[htbp]
	\centering
	\includegraphics[totalheight=3in]{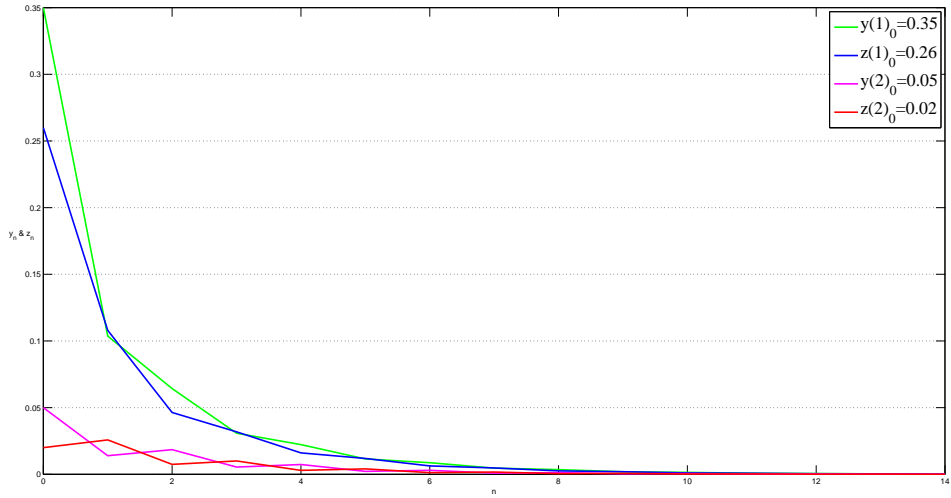}
	\caption{The comparation of solutions of system (4.2) from $(y(1)_0,z(1)_0)$ and $(y(2)_0,z(2)_0)$.} \label{fig:4}
\end{figure}

\section{Conclusions}
\setcounter{equation}{0}
\vskip4mm
\par
(i) System(1.6) is always persistent and bounded for $p,q \in (0, \infty), a,b \in (0,1) $.\\

(ii) System(1.6) always has a zero equilibrium; If $	a < p,~~ b < q$, and $ab< \min(p,q),$ then (0,0) is global asymptotic stable. \\

(iii) If $ab>pq$ , System(1.6)  has a unique positive equilibrium $\bar{y},\bar{z} \in (y_*,y^*)\times (z_*,z^*)$; Moreover, if $e^{y_*} >  a, ~~ e^{z_*} >  b$ and $	a z^* \le \bar{y}\left(p+z^*\right) e^{y_*},\  ~ b y^* \le \bar{z}\left(q+y^*\right) e^{z_*}$,then $(\bar{y},\bar{z})$ is global asymptotic stable.

\section*{Conflict of interest}
 The authors declare that they have no competing interests.

\section*{Acknowledgment}
This work was financially supported by Guizhou  Scientific and Technological Platform Talents ([2022]020-1), Scientific Research Foundation of Guizhou Provincial Department of Science and Technology([2020]1Y008, [2022]021, [2022]026), and Scientific Climbing Programme of Xiamen University of Technology (XPDKQ20021).

\section*{Reference}
\bigskip

\noindent[1] 	 M. Kot,		   Elements of Mathematical Ecology ,							 Cambridge University Press, New York, 2001.														\\
\noindent[2]	  R. Beverton and S. Holt,				 On the dynamics of exploited fish populations,							   Fisheries Investigations, Ser 2, , \textbf 19  (1957), 1--533.												\\
\noindent[3]	   MDL. Sen,		 The generalized Beverton-Holt equation and the control of populations,									    Applied Mathematical Modelling , \textbf 32  (2008), 2312-2328.												\\
\noindent[4]	 X. Ding, W. Li ,		   Stability and bifurcation of numerical discretization Nicholson blowflies equation with delay  ,												 Discrete Dyn. Nat. Soc , \textbf 2006   (2006), 1-12. Article ID 19413.									\\
\noindent[5]	 A.J. Nicholson ,		    An outline of the dynamics of animal populations ,								 Aust. J. Zool , \textbf 2   (1954), 9-65 .													\\
\noindent[6]	  E. C. Pielou,			    An Introduction to Mathematical Ecology ,											  John Wiley $\&$ Sons, New York, 1965.									\\
\noindent[7]	  E. C. Pielou,			    Population and Community Ecology, Gordon and Breach ,											   New York, 1974.									\\
\noindent[8]	 Camouzis, E. , Ladas, G. 					 Periodically forced Pielou's equation,						    J. Math. Anal. Appl,  \textbf 333(1)  (2007), 117-127.												\\
\noindent[9]	  Kulenovi$\acute c $, M.R.S. , Merino, 					 Stability analysis of Pielou's equation with period-two coefficient,									    J. Difference Equ. Appl,  \textbf 13(5)  (2007), 383-406.									\\
\noindent[10]	 Ishihara, Keigo; Nakata, Yukihiko,			   On a generalization of the global attractivity for a periodically forced Pielou's equation ,											 J. Difference Equ. Appl. \textbf 18(3)  (2012),375-396.									\\
\noindent[11]	 Zhao, Houyu  			   Analytic invariant curves for an iterative equation related to Pielou's equation  ,											 J. Difference Equ. Appl. \textbf 19  (2013), no. 7, 1082-1092.									\\
\noindent[12]	 Nyerges, Gabor,			    A note on a generalization of Pielou's equation ,											 J. Difference Equ. Appl. \textbf 14(5)  (2008), 563-565.									\\
\noindent[13]	  E. El-Metwally, E.A. Grove, G. Ladas, R. Levins, M. Radin,						    On the difference equation $x_ n+1 =\alpha +\beta x_ n-1 e^ -x_ n  $ ,									  Nonlinear Anal, \textbf  47   (2001), 4623-4634.								\\
\noindent[14]	   I. Ozturk, F. Bozkurt, S. Ozen, 					 On the difference equation $y_ n+1 =\frac \alpha +\beta e^ -y_ n    \gamma +y_ n-1  $,     Appl. Math. Comput,  \textbf 181  (2006), 1387-1393. 																		\\
\noindent[15]	  Papaschinopoulos, G.; Ellina, G.; Papadopoulos, K. B. ,						    Asymptotic behavior of the positive solutions of an exponential type system of difference equations ,									  Appl. Math. Comput, \textbf 245   (2014), 181-190.								\\
\noindent[16]	  Papaschinopoulos, G.; Radin, M. ; Schinas, C. J. ,						    Study of the asymptotic behavior of the solutions of three systems of difference equations of exponential form ,									  Appl. Math. Comput, \textbf  218   (2012), 5310-5318.								\\
\noindent[17]	  A. Matsumoto, F. Szidarovszky,					    Asymptotic behavior of a delay differential neoclassical growth model ,									  Sustainability, \textbf 5   (2013), 440-455.									\\
\noindent[18]	 X. Ding, R. Zhang ,			    On the difference equation$x_ n+1 =(\alpha x_ n  +\beta x_ n-1 )e^ -x_ n  $ ,											 Adv. Difference Equ , 2008 (2008),  http://dx.doi.org/10.1155/2008/876936. 7pages.									\\
\noindent[19]	  L. Shaikhet,		   Stability of equilibrium states for a stochastically perturbed Mosquito population equation  ,												 Dyn. Contin. Discrete Impuls. Syst. Ser. B Appl.Algorithms , \textbf  21 (2)   (2014), . 185-196.									\\
\noindent[20]	 T. Awerbuch, E. Camouzis, G. Ladas, R. Levins, E.A. Grove, M. Predescu ,								   A nonlinear system of difference equations , linking mosquitoes, habitats and community interventions  ,						 Commun. Appl. Nonlinear Anal , \textbf 15 (2)   (2008) 77-88.									\\
\noindent[21]	  B. Iricanin, S. Stevic,			    On two systems of difference equations ,							 Discrete Dyn. Nat. Soc , \textbf  4   (2010) . (Article ID 405121).													\\
\noindent[22]	  G. Papaschinopoulos, N. Fotiades, C.J. Schinas,			    On a system of difference equations including exponential terms ,							 J. Difference Equ. Appl , \textbf 20 (5-6)   (2014), 717-732.													\\
\noindent[23]	   MDL Sen, S. Alonso-Quesada,					 Control issues for the Beverton-Holt equation in ecology by locally monitoring the environment carrying capacity: Nonadaptive and adaptive cases,						    Applied Mathematics and Computation , \textbf 215  (2009), 2616-2633.   												\\
\noindent[24]	   M. Bohner, S. Streipert,					 Optimal harvesting policy for the Beverton-Holt model,								    Mathematical Biosciences and Engineering , \textbf 13  (2016), 673-695.										\\
\noindent[25]	  V. L. Kocic, G. Ladas,					   Global behavior of nonlinear difference equations of higher order with application ,										  Kluwer Academic Publishers, Dordrecht, 1993.  								\\
\noindent[26]	  M. R. S. Kulenonvic, G. Ladas,					    Dynamics of second order rational difference equations with open problems and conjectures ,										  Chapaman $\&$ Hall/CRC,Boca Raton, 2002.								\\
\noindent[27]	  Khan, A. Q.; Qureshi, M. N.,					    Behavior of an exponential system of difference equations ,								  Discrete Dyn. Nat. Soc, Art. ID 607281, 9 pp  (2014).										\\
\noindent[28]	  Papaschinopoulos, G.; Fotiades, N. ; Schinas, C. J.,						    On a system of difference equations including negative exponential terms ,									  J. Difference Equ. Appl., \textbf 20   (2014), 717-732.								\\
\noindent[29]	  W. Chen, W. Wang,					    Global exponential stability for a delay differential neoclassical growth model ,										  Adv. Difference Equ., \textbf 2014   (2014), 325.								\\
\noindent[30]	 L. Shaikhet ,		   Lyapunov Functionals and Stability of Stochastic Functional Differential Equations  ,												 Springer ,  Dordrecht, Heidelberg, New York, London,									\\
\noindent[31]	     Enatsu, Y, Nakata, Y, Muroya, Y,					      Global stability for a class of discrete SIR epidemic models ,									   Math. Biosci. Eng., \textbf 7 (2)   (2010), 347-361.									\\

\end{document}